\documentclass[12pt]{article}

\usepackage{amsmath}
\usepackage{amsfonts}
\usepackage{amssymb}
\usepackage{graphicx}

\newtheorem{proposition}{Proposition}[section]
\newtheorem{theorem}{Theorem}[section]
\newtheorem{lemma}[theorem]{Lemma}
\newtheorem{corollary}[theorem]{Corollary}
\newtheorem{remark}[theorem]{Remark}
\newtheorem{example}[theorem]{Example}

\newtheorem{definition}{Definition}

 \def\Xint#1{\mathchoice
       {\XXint\displaystyle\textstyle{#1}}%
       {\XXint\textstyle\scriptstyle{#1}}%
       {\XXint\scriptstyle\scriptscriptstyle{#1}}%
       {\XXint\scriptscriptstyle\scriptscriptstyle{#1}}%
       \!\int}
    \def\XXint#1#2#3{{\setbox0=\hbox{$#1{#2#3}{\int}$}
         \vcenter{\hbox{$#2#3$}}\kern-.5\wd0}}
    
    \def\dashint{\Xint-}

\def\phi{{\varphi}}

\DeclareSymbolFont{AMSb}{U}{msb}{m}{n}
\DeclareMathSymbol{\N}{\mathbin}{AMSb}{"4E}
\DeclareMathSymbol{\Z}{\mathbin}{AMSb}{"5A}
\DeclareMathSymbol{\R}{\mathbin}{AMSb}{"52}
\DeclareMathSymbol{\Q}{\mathbin}{AMSb}{"51}
\DeclareMathSymbol{\I}{\mathbin}{AMSb}{"49}
\DeclareMathSymbol{\C}{\mathbin}{AMSb}{"43}

\begin{document}

\title{Uniqueness of minimizers of weighted least gradient problems arising in conductivity imaging}

\author{{Amir Moradifam\footnote{Department of Mathematics, Michigan State University, East Lansing, MI, USA. E-mail: amir@math.msu.edu. }
\qquad Adrian Nachman\footnote{Department of Mathematics and the
Edward S. Rogers Sr. Department of Electrical and Computer
Engineering, University of Toronto, Toronto, Ontario, Canada M5S
2E4. E-mail: nachman@math.toronto.edu. }\qquad Alexandru
Tamasan\footnote{Department of Mathematics, University of Central
Florida, Orlando, FL, USA. E-mail: tamasan@math.ucf.edu. }}}
\date{}

\smallbreak \maketitle

\begin{abstract}
We prove uniqueness for minimizers of the weighted least gradient problem
\[\inf \left\lbrace \int_{\Omega} a|Du|: \ \ u\in BV(\Omega), \ \ u|_{\partial \Omega}=f  \right\rbrace.\]
The weight function $a$ is assumed to be continuous and it is allowed to vanish in certain subsets of $\Omega$. Existence is assumed a priori.  Our approach is motivated by the hybrid inverse problem of imaging electric conductivity from interior knowledge (obtainable by MRI) of the magnitude of one current density vector field.

\end{abstract}
\maketitle

\section{Introduction}
Consider the following  weighted least gradient problem
\begin{equation}\label{gen-prob}
\inf \left\lbrace \int_{\Omega} a|Du|: \ \ u\in BV(\Omega), \ \ u|_{\partial \Omega}=f  \right\rbrace,
\end{equation}
where $\Omega$ is a bounded open set in $\R^n$ ($n\geq 2$) with connected Lipschitz boundary, $a$ is a bounded non-negative function, and $f\in C(\partial \Omega)$.  
Our motivation comes from a hybrid inverse problem in medical imaging. The problem is to determine the conductivity of a body from knowledge of the magnitude $a=|J|$ (in $\Omega$) of one current density vector field $J$ generated by imposing the voltage $f$ on $\partial \Omega$. The interior data $|J|$ can be obtained non-invasively via a magnetic resonance technique pioneered in \cite{joy89}. In \cite{NTT08} this problem was reduced to the weighted least gradient problem (\ref{gen-prob}) (see \cite{NTT10, MNTa_SIAM, MNT, nashedTa10} for further results on partial data, inclusions, reconstruction algorithms, and stability, and also \cite{NTT11} for a review).

Existence and uniqueness of the minimizers of (\ref{gen-prob}) was first studied for the case $a\equiv 1$ in \cite{sternberg_ziemer92}  (see also \cite{sternbergZiemer93}). In particular  the authors proved that (\ref{gen-prob}) has a unique minimizer if $f$ is continuous and the mean curvature of $\partial \Omega$ is positive on a dense subset of $\partial \Omega$, see conditions (3.1) and (3.2) in \cite{sternberg_ziemer92}. Recently, in a companion paper \cite{JMN},  authors showed among other results that if $a\in C^{1,1}(\Omega)$ is positive and bounded away from zero, and if $f$ is continuous on $\partial \Omega$, then the weighted least gradient problem (\ref{gen-prob}) has at most one minimizer in $BV(\Omega)$. They also showed that the condition $a \in C^{1,1}(\Omega)$ is sharp in the sense that uniqueness can fail for $a \in C^{1,\alpha}(\Omega)$ with any $\alpha <1$.

In this paper the weight $a$ is only assumed to be continuous and it is allowed to vanish in certain subsets of $\Omega$. On the other hand here we require existence of a minimizer $u$ of (\ref{gen-prob}) that has appropriate properties (see Definition \ref{admis}). This assumption is naturally satisfied in the weighted least gradient problems arising in the conductivity problems (explained below) that motivated us. Our uniqueness proof is quite different from that in \cite{JMN}, and is based on a calibration argument. 

To motivate the existence assumption, assume $\Omega \subset \R^n$ is a conductive body with (spatially varying) conductivity $\sigma$. If the voltage $f$ is imposed on $\partial \Omega$, then the corresponding voltage potential $u$ is the solution of the following Dirichlet problem 
\begin{equation}\left\{ \begin{array}{ll}
\nabla\cdot \sigma \nabla u=0,&\mbox{in}\,\Omega,\\
u=f &\mbox{on}\,\partial \Omega.
\end{array}\label{pde_inclusions0} \right.
\end{equation}
Let $J=-\sigma \nabla u$ be the corresponding current density. In the inverse problem $\sigma$ is not known. It is shown in \cite{NTT08} that if $u$ satisfies (\ref{pde_inclusions0}) then it is a minimizer of the weighted least gradient problem (\ref{gen-prob}) which only involves the measured data $a=|J|$ and the prescribed voltage $f$ on $\partial \Omega$.

More generally, as in \cite{MNTa_SIAM}, we also consider  the case when $\Omega$ contains perfectly conducting and insulating inclusions $O_\infty$ and $O_0$.   In this case the corresponding voltage potential $u$ is the unique solution of the following equation 

\begin{equation}\left\{ \begin{array}{ll}
\nabla\cdot \sigma \nabla u=0,&\mbox{in}\,\Omega\setminus\overline{O_\infty \cup O_0},\\
\nabla u=0, &\mbox{in} \ \ O_\infty,\\
u|_+=u|_-,&\mbox{on}\ \ \partial (O_\infty \cup O_0),\\
\int_{\partial O_{\infty}^j}\sigma\frac{\partial u}{\partial \nu}|_{+}ds=0,&j=1,2,...,\\
\frac{\partial u}{\partial \nu}|_{+}=0,&\mbox{on}\;\partial O_0,\\
u|_{\partial\Omega}=f,
\end{array}\label{pde_inclusions} \right.
\end{equation}
where $O_0 \cap O_\infty=\emptyset$ and $O_\infty=\bigcup _{j=1} O_{\infty}^{j}$ is the partition of $O_\infty$ into open connected components (see the appendix in \cite{MNTa_SIAM} for more details).  Moreover, if $\sigma \in C^{\alpha}(\Omega \setminus \overline{O_0 \cup O_\infty})$ and the boundaries of $O_0$, $O_\infty$, and $\Omega$ are regular enough, then it follows from standard elliptic regularity results that $u\in C^{1}(\bar{\Omega} \setminus (O_0 \cup O_\infty))$.  Under certain assumptions, it is shown in \cite{MNTa_SIAM} that the solution of the equation (\ref{pde_inclusions}) is a minimizer of (\ref{gen-prob}), where $a$ is the magnitude of the corresponding  current density vector field. 

Uniqueness of minimizers in $W^{1,1}(\Omega)\cap C(\bar{\Omega})$ was proved in \cite{NTT08}, and \cite{MNTa_SIAM} in the presence of inclusions. The main objective of this paper is to prove uniqueness of minimizers of the above problem in $BV(\Omega)$ where we have compactness (see Proposition \ref{compact}). This is crucial when one studies the stability of the problem with respect to errors in measurements of $|J|$ and $f$. Once $u$ is determined, the shape and locations of perfectly conducting and insulating inclusions and the conductivity outside of the inclusions can be easily recovered.

We now state our assumptions and results more precisely. Throughout the paper we shall assume that $\Omega \subset \R^n$ is a bounded open set with connected Lipschitz boundary and $a\in L^{\infty}(\Omega)$ is non-negative. By $\mathcal{H}^{d}$ we denote the $d$-dimensional Lebesgue/Hausdorff measure. While $a$ is allowed to vanish, its zero set
\[S:=\{x\in \bar{\Omega}: a(x)=0 \}\] is assumed to satisfy the following {\em structural hypothesis}:
\begin{equation}\label{str.hyp}
\bar{S}=:\bar{O}_0 \cup \Gamma,
\end{equation}
where $\Gamma$ is a set of measure zero and the (possibly empty) open set $O_0 \subset \subset \Omega$, modelling the insulating regions, is a pairwise disjoint union of finitely many $C^1$- diffeomorphic images of the unit ball. In addition, in two dimensions $O_0$ is assumed to have at most one such component.

Let
\[ X:= \{b\in L^{\infty}(\Omega;\mathbb{R}^n): \; \nabla \cdot b \in L^n(\Omega)\}.\]
For any $u \in BV_{loc}(\Omega \setminus \bar{S})$ the total variation of $u$ (with respect to $a$)
in $\Omega$ is defined as
\begin{equation}\label{BV-def}
\int_{\Omega} |Du|_{a}=\sup_{b \in \mathfrak{B}_a} \int_{\Omega} u  \nabla \cdot b  \ \ dx,
\end{equation}
where
\begin{equation*}
\mathfrak{B}_a=\{ b\in X: \hspace{0.15cm} \mbox{supp}(b)\subset \subset \Omega,\; |b(x)|\leq a(x)\; \mathcal{H}^{n}\hbox{-a.e. in} \ \ \Omega\},
\end{equation*}
(see \cite{AB} and the references cited therein). By the structural hypothesis (\ref{str.hyp}), $\partial S$ has measure zero and therefore $\int_{\Omega} |D u|_a$ is independent of the value of $u$ in $\bar{S}$. Hence $BV_{loc}(\Omega \setminus \bar{S})$ is the natural space of functions in which (\ref{BV-def}) makes sense.\\

Now consider the weighted least gradient problem
\begin{equation}\label{mainProb}
\min \{\int_{\Omega}|Dv|_a:  v\in BV_{loc}(\Omega\setminus \bar{S}),\; v|_{\partial\Omega}= f \},
\end{equation}
where the boundary condition is in the sense of the trace of functions in $BV(\Omega)$. In general the minimization problem $(\ref{mainProb})$ need not have a unique solution (see \cite{JMN}). The following admissibility assumption plays a crucial role in our uniqueness proof. Essentially it assumes continuity of $a$ outside of the inclusions and existence of a minimizer of (\ref{mainProb}).

\begin{definition}[Admissibility]\label{admis}
 Let $\Omega \subset \R^n$ be a open bounded region with connected Lipschitz boundary. A pair of functions $(f,a)$ is called {\em admissible} if the following conditions hold.

 (i) The zero set $S$ of the weight $a$ satisfies the structural hypothesis \eqref{str.hyp} for some $O_0$ and $\Gamma$.

 (ii) There exists a solution $u\in C^1(\bar{\Omega}\setminus O_0)$  of the weighted least gradient problem (\ref{mainProb}) such that
 \[\{x\in \overline{\Omega}: |\nabla u (x)|=0 \} \setminus \bar{O}_0=O_{\infty},\] 
 and  int$(\overline{ u(O_\infty)})=\emptyset$. 

 (iii) $a \in  C(\bar{\Omega} \setminus (O_0 \cup O_{\infty}))$.
 \end{definition}

The (possibly empty) set $O_\infty$ models the perfectly conducting inclusions. Note that the above definition of admissibility is significantly simplified if $O_0=O_\infty=\emptyset$. Even in this simpler case  a large class of admissible pairs $(f,a)$ is provided by the conductivity problem (\ref{pde_inclusions0}). 

\medskip

The following is our uniqueness result.

\medskip

\begin{theorem}[Uniqueness]\label{mainResult}
Let $\Omega \subset \R^n$ be a bounded Lipschitz domain with connected boundary and $(f,a)$ be admissible in the sense of Definition \ref{admis}. Then the weighted least gradient problem (\ref{mainProb}) has a unique solution in $BV_{loc}(\Omega \setminus \bar{S})$.
\end{theorem}

It may be helpful to compare the above theorem to the uniqueness result in \cite{JMN}. The uniqueness proof in \cite{JMN} does not require the pair $(a,f)$ to be admissible, but it only works for $0<c<a\in C^{1,1}(\Omega)$.\\

To illustrate a simple case with one perfectly conducting inclusion, consider the following example from \cite{sternberg_ziemer}. 
\begin{example}\label{ex}
Let  $D=\{x\in R^2: \ \ x^2+y^2<1\}$ be the unit disk, $f(x,y)=x^2-y^2$, and $O_\infty=(-\frac{1}{\sqrt{2}}, \frac{1}{\sqrt{2}} )\times (-\frac{1}{\sqrt{2}}, \frac{1}{\sqrt{2}})$.  It is shown in \cite{sternberg_ziemer} (see also \cite{MNTa_SIAM} for a different proof) that 
\begin{eqnarray*}
u=\left\{ \begin{array}{ll}
2x^2-1, \ \ \text{if } \ \ |x|\geq \frac{1}{\sqrt{2}}, \ \ |y|\leq \frac{1}{\sqrt{2}},\\
0, \hspace{0.5in}\ \hbox{if} \ \ (x,y)\in O_\infty,\\
1-2y^2, \ \ \text{if } \ \ |x|\leq \frac{1}{\sqrt{2}}, \ \ |y|\geq \frac{1}{\sqrt{2}}.\\
\end{array} \right.
\end{eqnarray*}
is a minimizer of the least gradient problem 
\begin{equation}\label{sz}
\hbox{min}\{\int_{D}|\nabla u| dx, \ \ u\in BV(D), \ \ \hbox{and}\ \ u|_{\partial D}=f\}.
\end{equation}
It is easy to observe that $(1,x^2-y^2)$ is an admissible pair with $O_\infty$ defined as above and $S=O_0=\emptyset$.  Hence Theorem \ref{mainResult} provides a new proof that $u$ is the only minimizer in $BV(\Omega)$.  
\end{example}

\section{Preliminaries}
In this section we recall and present some preliminary results that will be used in the following sections. First we recall a useful representation formula from \cite{AB}.  For  $u\in BV(\Omega)$
\begin{equation}\label{BV-defAB}
\int_{A}|Du|_a=\int_{A} h(x,v^u)|Du|,
\end{equation}
where
\begin{equation}
h(x,v^u)=(|Du|-\hbox{ess} \sup_{b\in \mathfrak{B_a}}) b \cdot v^u(x) \ \ \ \ for \ \ |Du|-a.e. \ \ x\in \Omega
\end{equation}
and $v^u$ denotes the Radon-Nikodym derivative $v^{u}(x)=\frac{d\, Du}{d\, |Du|}$.
The right-hand side of \eqref{BV-defAB} makes sense, since $v^u$ is $|Du|$-measurable,
and hence $h(x, v^u(x))$ is as well. In particular, if $u\in BV(\Omega)$, and the coefficient $a$ is continuous in the Borel measurable subset $A \subset \Omega$, then
\begin{equation}
\int_{A}|Du|_a=\int_{A} a|Du|,
\end{equation}
as shown in \cite{AB}.  The following Lemma provides a simple extension of this formula for the total variation of the voltage potential $u$ that corresponds to an admissible pair $(f,a)$.

\begin{lemma}\label{totalVaVP}
 Let $\Omega \subset \R^n$ be a bounded open region with Lipschitz boundary, $(f,a)$ be admissible, and $u$ be a minimizer of (\ref{mainProb}) as in Definition \ref{admis}. Then
\[\int_{\Omega}|D u|_a=\int_{\Omega}a |\nabla u| dx. \]
\end{lemma}
{\bf Proof:}  Since $(f,a)$ is admissible, $a \in C^0(\Omega \setminus (O_0 \cup O_{\infty}))$.  Hence by \cite[Proposition 7.1]{AB} we have that
\begin{equation}\label{hFormula}
h(x,v^u)= \left\{ \begin{array}{ll}
a(x) \ \ &\hbox{in} \ \  \Omega \setminus \overline{(O_0\cup O_{\infty})}\\
0 \ \  &\hbox{in}\ \ O_0. \\
\end{array} \right.
\end{equation}
Thus it follows from $(\ref{BV-defAB})$ that
\begin{eqnarray*}
\int_{\Omega} |Du|_a=\int_{\Omega \setminus (O_0\cup O_{\infty})} a|\nabla u|=\int_{\Omega} a|\nabla u|dx.
\end{eqnarray*}
\hfill $\Box$

It is a straightforward consequence of the definition (\ref{BV-def}) that $u\mapsto\int_{\Omega} |Du|_{a}$ is $L^{\frac{n}{n-1}}(\Omega)-$lower semi-continuous. As proven in \cite[Theorem 1.2]{An},  if $\nu$ denotes the outer unit normal vector to $\partial \Omega$, then for every $b\in X$ there exists a unique function $[b\cdot \nu] \in L^{\infty}_{\mathcal{H}^{n-1}}(\partial \Omega)$ such that
\begin{equation}
\int_{\partial \Omega} [b\cdot \nu] u d\mathcal{H}^{n-1}=\int_{\Omega} u \nabla \cdot b dx +\int_{\Omega}  b \cdot \nabla u dx, \ \ \forall u \in C^1(\bar{\Omega}).
\end{equation}

Moreover, for $u\in BV(\Omega)$ and $b\in X$, the linear functional $u\mapsto(b \cdot Du)$ gives rise to a Radon measure on $\Omega$, and
\begin{equation}\label{IBP0}
\int_{\partial \Omega} [b\cdot \nu]u d\mathcal{H}^{n-1}=\int_{\Omega} u \nabla \cdot b dx +\int_{\Omega}  (b \cdot D u), \ \ \forall u \in BV(\Omega),
\end{equation}
see \cite{Al, An} for a proof. We shall need the following lemma in the proof of our uniqueness result.

\begin{lemma}\label{IBP} Let $S$ be as defined in (\ref{str.hyp}) and $b\in X$. If $u\in L^{\infty}(\Omega)$ and $\int_{\Omega}|Du|_a<\infty$, then 
\begin{equation}\label{IBP-eq}
\int_{\partial \Omega} [{b}\cdot v_{\Omega}]u d\mathcal{H}^{n-1}=\int_{\Omega} u \nabla \cdot {b} dx +\int_{\Omega}  ({b} \cdot D u),
\end{equation}
for some unique function $[b\cdot \nu]$ in $L^{\infty}_{\mathcal{H}^{n-1}}(\partial \Omega)$.  
\end{lemma}
{\bf Proof:} By the structural hypothesis (\ref{str.hyp}), $S$ has finite perimeter in $\Omega$. Therefore, it follows from $\int_{\Omega}|Du|_a<\infty$ that 
\[
BV_{loc}(\Omega \setminus \bar{S}) \cap L^{\infty}(\Omega) \subset BV(\Omega). 
\]
Now (\ref{IBP-eq}) follows from (\ref{IBP0}). \hfill $\Box$
 
\medskip 
 
The following compactness result shows that $BV_{loc}(\Omega \setminus \bar{S})$ is the natural space of function for the minimization problem (\ref{mainProb}). 
\begin{proposition}[Compactness] \label{compact} Let $a\in L^{\infty}(\Omega) $ and assume that the set
\[S=\{x\in \Omega: a(x)=0 \} \]
satisfies the structural hypothesis (\ref{str.hyp}). Then every sequence $\{u_n\}_{n=1}^{\infty}$ in $L^1(\Omega \setminus \bar{S})$ with
\[c:=\sup_n\int_{\Omega} |Du_n|_a<\infty\]
has a convergent subsequence in $L^1_{loc}(\Omega \setminus \bar{S})$. Moreover if $u$ is a limit point, then $u\in BV_{loc}(\Omega\setminus \bar{S})$.
\end{proposition}
{\bf Proof:}
Consider the nested exhaustion of $\Omega\setminus \bar{S}$ by the
open subsets
\begin{align}\label{exhaustion}
\Omega^{k}:=\{x\in\Omega:\;a(x)>1/k\}, \quad k\in\N,
\end{align}i.e. $\Omega^{k}\subset\Omega^{k+1}$ and
$ \cup_{k=1}^\infty\Omega^{k}=\Omega\setminus \bar{S}$. For each fixed $k\in N$
\[\int_{\Omega^{k}}|Du_n| \leq k \int_{\Omega^{k}} a |D u|dx =\int_{\Omega^k}|Du|_a\leq kc, \ \ \hbox{for all} \ \ n\in N. \]

The classical compactness embedding of $BV(\Omega^{k})$ in
$L^1(\Omega^{k})$ yields a subsequence
$\{u_{n^{1}_{i}}\}_{i=1}^{\infty}$ of $\{u_i\}_{i=1}^{\infty}$ such
that $u_{n^{1}_{i}} \rightarrow v_1$ in $L^1(\Omega^1)$. Similarly,
there exists a subsequence $\{u_{n^{2}_{i}}\}_{i=1}^{\infty}$ of
$\{u_{n^{1}_{i}}\}_{i=1}^{\infty}$ such that $u_{n^{2}_{i}}
\rightarrow v_2$ in $L^1(\Omega^{1/2})$, and $v_1= v_2$ on
$\Omega^1$.  Repeating this argument we obtain a family of subsequences (indexed
in $k$) $\{u_{n^{k}_{i}}\}_{i=1}^{\infty}$ such that
$u_{n^{k}_{i}} \rightarrow v_k$ in $L^1(\Omega^{k})$, for each $k$
fixed. Since $\cup_{k=1}^\infty\Omega^{k}=\Omega\setminus \bar{S}$ and
$v_j= v_k$ on $\Omega^{k}$ for all $j\geq k$, one can define
a function $u$ on $\Omega\setminus \bar{S}$ by setting $u:=v_k$ in each
$\Omega^{k}$. Any compact $K \subset \Omega\setminus \bar{S}$ is
contained in $\Omega^{k}$ for $k$ large enough, hence
$\{u_{n^{k}_{i}}\}_{i=1}^{\infty}$ converges to $u$ in $L^1(K)$.  Since $\int_{\Omega}|Du|_a$ is lower semi-continuous,  $\int_{\Omega} |Du|_a \leq c.$ \hfill $\Box$\\

The next two results yield a calibration which will be used in the uniqueness proof. \\
Suppose $a\in L^2(\Omega)$ and fix $u_f \in H^1(\Omega)$ with $u_f|_{\partial \Omega}=f$. Consider the weighted least gradient problem

\begin{equation*}
(P) \ \ \ \ \min_{v\in H^1_0(\Omega)} \int_{\Omega} a|\nabla v+\nabla u_f|dx.
\end{equation*}

In \cite{MNT} it is shown that the dual problem to $(P)$ is  
\begin{equation*}
(D) \ \ \ \ \max \{<\nabla u_f, b>: \ \ b\in L^2(\Omega;\mathbb{R}^n),\ \ |b(x)|\leq a(x) \ \ a.e. \ \ \hbox{and} \ \ \nabla \cdot b\equiv 0\}.
\end{equation*}
Let $v(P)$ and $v(D)$ be the optimal values of the primal and dual problems. It is shown in \cite{MNT} that $v(P)=v(D)$ and the dual problem $(D)$ has an optimal solution. The following proposition is an immediate  consequence of Proposition 2.1 and Corollary 2.3 in \cite{MNT}.

\begin{proposition}\label{DualTheo}
Let $a\in L^2(\Omega)$ be a non-negative function and $v_f \in H^1(\Omega)$ with $v_f|_{\partial \Omega}=f$. Then the optimal values of the primal problem $(P)$ and dual problem $(D)$ are equal, and the dual problem $(D)$ has an optimal solution $J$ with $\nabla \cdot J\equiv 0$ in $\Omega$. Moreover, if $v$ is an optimal solution of the primal problem $(P)$, then
\[J(x)=a(x)\frac{\nabla( v(x)+v_f(x))}{|\nabla (v(x)+ v_f(x))|} \ \ \hbox{if} \ \ |\nabla( v(x)+v_f(x))|\neq 0,\]
for all $x\in \Omega$.
\end{proposition}

The following result is an immediate consequence of Proposition \ref{DualTheo}.

\begin{corollary}\label{PropDual}
Let $\Omega \subset \R^n$ be a bounded Lipschitz domain and $(f,a)$ be an admissible pair. Then there exists an optimal solution $J\in L^{2}(\Omega;\mathbb{R}^n)$ of the dual problem (D) such that $\nabla \cdot J\equiv 0$ in $\Omega$, $|J|\leq a$ a.e. in $\Omega$, and with $u$, $O_0$, and $O_{\infty}$ as described in Definition \ref{admis} we have

\begin{equation}
J(x)= \left\{ \begin{array}{ll}
a(x)\frac{\nabla u}{|\nabla u|}\ \ &\hbox{if} \ \  |\nabla u|\neq 0\\
0 \ \  &\hbox{if}\ \ a(x)=0. \\
\end{array} \right.
\end{equation}
Moreover $J$ is continuous in $\bar{\Omega} \setminus \overline{(O_0 \cup O_\infty)}$, and $|J(x)|>0$ whenever $a(x)>0$.

\end{corollary}

\section{Uniqueness of minimizers}

In this section we prove our main result, Theorem \ref{mainResult}. Let $u$ be the minimizer of the weighted least gradient problem (\ref{mainProb}) assumed in the Definition \ref{admis}, and  suppose $u_1\in BV_{loc}(\Omega \setminus \bar{S})$ is another minimizer.  We will show that $u=u_1$ a.e. in $\Omega \setminus \bar{S}$. First notice that $u_1$ is bounded above and below almost everywhere. Indeed if we define

\begin{equation}\label{mM}
\bar{u}_1(x)= \left\{ \begin{array}{ll}
u(x)\ \ &\hbox{if} \ \  m_f\leq u_1(x) \leq M_f\\
M_f \ \  &\hbox{if}\ \ u_1(x)>M_f, \\
m_f\ \ &\hbox{if}\ \ u_1(x)<m_f,
\end{array} \right.
\end{equation}
where $M_f$ and $m_f$ are the maximum and minimum values of $f$ on $\partial \Omega$, then it is easy to see that $\bar{u}_1 \in BV_{loc}(\Omega\setminus \bar{S})$ and
\begin{equation}
\int_{\Omega}|D\bar{u}_1|_a \leq \int_{\Omega}|Du|_a.
\end{equation}
Moreover the inequality is strict if $\{x\in \Omega: u_1(x)>M_f\}$ or $\{x\in \Omega: u_1(x)<m_f\}$ has positive measure.  Therefore we may assume $range(u_1)\subset range(f)$.\\

Next we prove that
\[ \frac{\nabla u}{|\nabla u|}=\frac{d Du_1}{d|Du_1|}\]
$|Du_1|-a.e.$ in $\Omega \setminus \overline{O_0 \cup O_{\infty}}$.

\begin{lemma}
Let $(f,a)$ be an admissible pair and $u$ be the corresponding minimizer of $(\ref{mainProb})$. If $u_1$ is another minimizer, then
\[\frac{\nabla u}{|\nabla u|}=v^{u_1} \ \ \ \ |Du_1|-a.e.\ \ \hbox{in} \ \  \Omega \setminus \overline{O_0 \cup O_{\infty}}.\]
\end{lemma}
{\bf Proof:} Let $x\in \Omega$ and choose $\epsilon>0$ small enough such that $B(x, 2\epsilon) \subset \Omega$. Then it follows from the definition of $h(x, v^{u_1})$ that
\[\int_{B(x, \epsilon)} h(x,v^{u_1})|Du_1|\geq \int_{B(x, \epsilon)} J \cdot v^{u_1}|Du_1|,\]
where $J$ is the solution of the dual problem $(D)$ in Proposition \ref{PropDual}. Therefore
\[h(x,v^{u_1})\geq J\cdot v^{u_1}, \ \ \ \ |Du_1|-a.e. \ \ \hbox{in} \ \ \Omega.\]
Thus
\begin{eqnarray*}
\int_{\Omega} |Du_1|_a&=&\int_{\Omega}h(x, v^{u_1})|Du_1| \geq \int_{\Omega} J \cdot v^{u_1} |Du_1|\\
&=& \int_{\Omega} J\cdot Du_1= \int_{\partial \Omega} J\cdot \nu f d\mathcal{H}^{n-1}\\
&=&\int_{\Omega}\nabla u\cdot J dx= \int_{\Omega}|J| |\nabla u| \\
&=& \int_{\Omega} |Du|_a=\int_{\Omega} |Du_1|_a,
\end{eqnarray*}
where the third and fifth equalities follow form Lemma \ref{IBP} and Lemma \ref{totalVaVP}, respectively. Therefore
 \[h(x,v^{u_1})=J\cdot v^{u_1},\ \ |Du_1|-a.e. \ \ \hbox{in} \ \ \Omega.\]

Since $a$ is continuous in $\Omega \setminus \overline{(O_0 \cup O_{\infty})}$, as in (\ref{hFormula}) we have
\[ h(x, v^{u_1})=a(x), \ \  |Du_1|-a.e. \ \ \hbox{in} \ \  \Omega \setminus \overline{(O_0 \cup O_{\infty})}.\]
On the other hand $|v^{u_1}|=1$ and $|J|\leq a$, $|Du_1|-a.e.$ in $\Omega$, and $|\nabla u|\neq 0$ on $\Omega \setminus \overline{O_0 \cup O_{\infty}}$. Thus
\[\frac{\nabla u}{|\nabla u|}=\frac{J}{|J|} =v^{u_1}, \ \ |Du_1|-a.e. \ \ \hbox{in} \ \ \Omega \setminus \overline{O_0 \cup O_{\infty}}.\]
\hfill $\Box$

For $\lambda \in range(u_1)$, let
\[E_{\lambda}=\{x\in \Omega \setminus \overline{O}_0: \ \ u_1(x) \geq t\}.\]
Define
\begin{equation}\label{Eprime}
E'_\lambda:=\{x\in \R^n: \lim_{r \rightarrow 0} \frac{\mathcal{H}(B(r,x)\cap E_\lambda)}{\mathcal{H}(B(r))}=1\}.\
\end{equation}
By changing $u_1$ in a set of measure zero, we may assume that $E_\lambda=E'_{\lambda}$. Indeed throughout this paper we shall always assume that $E_\lambda=E'_{\lambda}$. We also define 
\begin{equation}\label{Z}
Z=\{x\in \bar{\Omega}\setminus O_0: u(x)\in \overline{u(O_\infty)}\},
\end{equation}
where $u$ is the minimizer of (\ref{mainProb}) in Definition \ref{admis}.

\begin{lemma}\label{sameLevelSets}
Assume that $(f,a)$ is an admissible pair, $u$ is the corresponding minimizer of $(\ref{mainProb})$, and $u_1 \in BV_{loc}(\Omega \setminus \bar{S})$ is another minimizer. Let $\Sigma$ be a connected component of $E_\lambda$, then for almost every $\lambda \in range (u_1)$, either \\
(i) $\Sigma \subset Z$\\
or \\
(ii)  $\Sigma \cap Z=\emptyset$, $\Sigma$ is a $C^1$ hypersurface, and  $u$ is constant on  $\Sigma$.\\
\end{lemma}

Let $\Lambda$ be the set of all $\lambda \in range(u_1)$ such that every connected component $\Sigma$ of $E_\lambda$ with $\Sigma \cap Z=\emptyset$ is a $C^1$ hypersurface. Then by the above lemma 
\begin{equation}\label{Lambda}
\mathcal{H}^1(range(u_1)\setminus \Lambda)=0.
\end{equation}

{\bf Proof of Lemma \ref{sameLevelSets}:} By co-area formula we have
\begin{align}
0=\int_{\Omega\setminus \overline{O_0 \cup O_\infty}} \varphi [\frac{\nabla u}{|\nabla u|} - v^{u_1}] |Du_1| = \int_{0}^{\infty} \int_{\partial^* E_{\lambda} \cap (\Omega\setminus \overline{O_0 \cup O_\infty})} \varphi [\frac{\nabla u}{|\nabla u|} -v^{u_1}]d \mathcal{H}^{n-1} d\lambda
\end{align}
for every smooth vector field $\varphi$, where $\partial^*E_\lambda$ is the reduced boundary of $E_\lambda$. Therefore $\frac{\nabla u}{|\nabla u|}=v^{u_1}$, $\mathcal{H}^{n-1}-a.e.$ in $\partial^* E_\lambda \cap (\Omega \setminus \overline{O_0 \cup O_\infty})$ for almost every $\lambda \in range (u_1)$.  Since $|D\chi_{E_{\lambda}}|$ is the $(n-1)-$dimensional Hausdorff measure restricted to $\partial^* E_\lambda$ (see \cite{G}, Chapter 4), for almost every $\lambda \in range(u_1)$,  the generalized normal $\nu(x)$ exists for $|D\chi_{E_{\lambda}}|-a.e.$  $x\in \partial E_\lambda \cap (\Omega \setminus \overline{O_0 \cup O_\infty})$ and coincides there with the continuous vector field $\frac{\nabla u}{|\nabla u|}$ . By Theorem 4.8 in \cite{G}, for every $x\in \partial E_\lambda \cap \Omega \setminus \overline{O_0 \cup O_\infty}$, $\partial E_\lambda$ can be represented as the graph of a Lipschitz continuous function $g$. Thus the derivative of $g$ coincides almost everywhere with a  continuous function and therefore $g$ must be $C^1$ and consequently we conclude that each connected component of $\partial E_\lambda \cap (\Omega \setminus \overline{O_0 \cup O_\infty})$ is a $C^1$ hyperspace for almost every $\lambda \in range(u_1)$.

Now we show that $u$ is constant on every $C^1$ connected component $\Sigma$
of $\partial E_\lambda \cap (\Omega \backslash \overline{O_0 \cup O_\infty}))$. Let
$\gamma:(-\epsilon,+\epsilon) \rightarrow
\Sigma$ be an arbitrary $C^1$ curve. Then

\[\frac{d}{dt}u(\gamma(s))=|\nabla u(\gamma(s))|\nu(\gamma (s)).\gamma'(s)=0,\]
because either $|\nabla u(\gamma(s))|=0$ or $\nu(\gamma
(s)).\gamma'(s)=0$ on $\Sigma$. Thus $u$ is constant along $\gamma$ and consequently $u$ is constant on $\Sigma$. The proof is now complete.
\hfill $\Box$\\

We show next that every connected component of $\partial E_{\lambda}$ intersects the boundary $\partial \Omega$.
\begin{proposition} \label{intersectBoundary}
Let $(f,a)$ be an admissible pair and $u_1$ be a minimizer of $(\ref{mainProb})$. Assume $\Sigma_{\lambda}$ is a $C^1$ connected component  of $\partial E_{\lambda}=\partial\{x \in \Omega\setminus O_0: \ \ u_1(x)>\lambda\}$, and $\Sigma_{\lambda}\cap Z=\emptyset$. Then
\[\overline{\Sigma}_{\lambda}\cap \partial \Omega\neq \emptyset.\]
\end{proposition}
{\bf Proof:} Assume $\overline{\Sigma}_{\lambda}\cap \partial \Omega = \emptyset$. We consider two cases:\\

(I) $\overline{\Sigma}_\lambda$ is a manifold without boundary in $\overline{\Omega} \setminus O_0$,

(II) $\overline{\Sigma}_\lambda \cap \partial O_0 \neq \emptyset$.\\

Case I: Assume that $\overline{\Sigma}_\lambda $ is a manifold without boundary in $\overline{\Omega}$. Then
$\partial \Omega \cup \Sigma_\lambda$ is a compact manifold with two
connected components. By the Alexander duality theorem for $\partial \Omega \cup \Sigma_\lambda$ (see, e.g., Theorem
27.10 in \cite{G81}) we have that $\R^n \setminus(\partial \Omega
\cup \Sigma_\lambda)$ is partitioned into three open connected components:
$\R^n=(\R^n\setminus \overline{\Omega}) \cup U_1\cup U_2$. Since
$\Sigma_\lambda \subset \Omega$ we have $U_1 \cup U_2=\Omega \setminus
\Sigma_\lambda$ and then $\partial U_i \subset
\partial \Omega \cup \Sigma_\lambda$ for $i=1,2$.

We claim that at least one of the $\partial U_1$ or $\partial U_2$
is in $\Sigma_\lambda$. Assume not, i.e. for $i=1,2$, $\partial U_i \cap
\partial \Omega\neq \emptyset$. Since $\partial \Omega$ is connected
(by assumption) we have that $U_1 \cup U_2\cup \partial \Omega$ is
connected which implies that $U_1 \cup U_2\cup(\R^n\setminus
\Omega)$ is also connected. Again by applying the Alexander duality
theorem for $\Sigma_\lambda \subset \R^n$, we have that $\R^n \setminus
\Sigma_\lambda$ has exactly two open connected components, one of which is
unbounded: $\R^n \setminus \Sigma_\lambda=U_{\infty}\cup U_0$. Since
$U_1\cup U_2 \cup (\R^n \backslash \Omega)$ is connected and
unbounded, we have that $U_1\cup U_2 \cup (\R^n \backslash \Omega)
\subset U_{\infty}$, which leaves $U_0 \subset \R^n \setminus (U_1
\cup U_2 \cup (\R^n \setminus \Omega)) \subset \Sigma_\lambda$. This is
impossible since $U_0$ is open and $\Sigma_\lambda$ is a hypersurface.
Therefore either $\partial U_1$ or $\partial U_2$ or both lie in
$\Sigma_\lambda$.

Assume $\partial U_1\subset \Sigma_t$. We claim that $u$ is constant in
$U_1$. Indeed, by Lemma \ref{sameLevelSets}, $u=c$ on $\Sigma_{\lambda}$ for some $c$. Hence the new map $\tilde{u}$
defined by

\begin{eqnarray*}
\tilde{u}:=\left\{ \begin{array}{ll}
u,\ \ x \in \Omega \setminus U_1,\\
c, \ \ \ \ x\in \overline{U_1},
\end{array} \right.
\end{eqnarray*}
is in $BV_{loc}(\Omega \setminus \bar{S})$ and decreases the
energy, which contradicts the minimality of $u$. Therefore
$u=c$ in $U_1$. This is a
contradiction since we have assumed $\overline{\Sigma}_{\lambda} \cap Z=\emptyset$.\\

Case II: Suppose $\overline{\Sigma}_\lambda \cap \partial O_0 \neq \emptyset$.  First assume $n\geq 3$. Let 
\[\epsilon^*:= min\{\min_{i\neq j} dist(O_i, O_j), \min_{i} dist(O_i, \partial \Omega)\},\]
where $O_i$, $1\leq i\leq m$, are the open connected components of the set $O_0$. For any  $0<\epsilon <\epsilon^*:$ define 

\[O_0^{\epsilon}=O_0\cup \{x\in \Omega: dis(x, O)<\epsilon\}.\]

Then $O_0^{\epsilon}$ is an open set with the same number of disjoint open connected components as $O_0$.  Now let $\Sigma_{\lambda}^{\epsilon}= \Sigma_\lambda \setminus O_0^{\epsilon}$ which we know is $C^1$ on $\Omega \setminus O_0^{\epsilon}$. Since $\partial \Sigma_{\lambda}^{\epsilon} \subset \partial O_{0}^{\epsilon}$ and  $\partial O_0^\epsilon \setminus \partial \Sigma_{\lambda}^{\epsilon} $ is open, each connected component of $\partial \Sigma_\lambda^{\epsilon}$ is the boundary of an open set in $\partial O_0^{\epsilon}$ with connected boundary.  Suppose $M$ is a connected component of $\partial \Sigma_{\lambda}^{\epsilon}$.  Then $M \subset \partial O^{\epsilon}_i$ for some $1 \leq i \leq m$, $O^{\epsilon}_i$ is $C^1$-diffeomorphic image of the unit ball for $\epsilon$ small, and $M$ is an orientable manifold without boundary in $\partial O_{0}^{\epsilon}$. Therefore it follows from Alexander's duality theorem that
\[ \partial O_i^{\epsilon} \setminus M= V_1 \cup V_2,\]
where $V_1, V_2$ are disjoint open connected (with respect to the topology of $\partial O_0^{\epsilon}$) sets.  Since $\Sigma^t_{\epsilon}$ can be extended to a $C^1$ hypersurface $\Sigma_\lambda$ inside $ O_0^{\epsilon} \setminus O_0$,  we can extend  $\Sigma_\lambda^{\epsilon}$ inside $O^{\epsilon}_i$ to obtain a $C^1$ hypersurface  $\Sigma$ such that
\[\Sigma \cap (\Omega \setminus O_0^{\epsilon})=\Sigma_\lambda^{\epsilon} \cap (\Omega \setminus O_0^{\epsilon})\]
and $\partial (\Sigma \cap O_0^{\epsilon})=M$. Repeating this argument for other connected components of $\partial \Sigma_\lambda^{\epsilon}$ leads to a $C^1$ orientable hypersurface $S^{\epsilon}$ with no boundary,
$\partial \Omega \cap S^{\epsilon}=\emptyset$, and $S^{\epsilon} \cap (\Omega \setminus O_0^{\epsilon}) =\partial \Sigma^{\epsilon}_\lambda$. Now apply Alexander's duality theorem to get the partition
\[\R^n \setminus S^{\epsilon}=U^{\epsilon} \cup U^{\epsilon}_{\infty},\]
where $U^{\epsilon}$ and $U^{\epsilon}_{\infty}$ are open subsets of $\R^n$ and  $U^{\epsilon}_{\infty}$ is unbounded. Notice that $\Sigma^{\epsilon}_\lambda \subset \partial U^{\epsilon} \subset \Sigma^{\epsilon}_\lambda \cup O_0^{\epsilon}$ and consequently $\Sigma^{\epsilon}_\lambda \subset \partial (U^{\epsilon} \setminus \bar{O}_0^{\epsilon}) \subset \partial O_0^{\epsilon} \cup \Sigma^{\epsilon}_t$.  If $\epsilon'< \epsilon$, then $\Sigma^{\epsilon}_\lambda \subset \Sigma^{\epsilon'}_\lambda$ and $O_0^{\epsilon '} \subset O_0^{\epsilon}$. Therefore
\[U^{\epsilon} \setminus \bar{O}_0^{\epsilon} \subset U^{\epsilon'} \setminus \bar{O}_0^{\epsilon'}.\]
Now let
\[U=\cup_{0<\epsilon <\epsilon^{*}} (U^{\epsilon} \setminus \bar{O}_0^{\epsilon}) . \]
Then $U$ is open and $\partial U \subset \Sigma_\lambda \cup O_0$.  We claim that $u$ is constant in $U$.  Indeed, by Lemma \ref{sameLevelSets} $u=c$ on $\Sigma_{\lambda}$ for some $c$ and the new map
defined by
\begin{eqnarray}\label{cont1}
\tilde{u}:=\left\{ \begin{array}{ll}
u,\ \ x \in \Omega \backslash  U,\\
\lambda, \ \ \ \ x\in U,
\end{array} \right.
\end{eqnarray}
is in $BV_{loc}(\Omega \setminus S)$. This contradicts the
minimality of $u$. Thus $u=c$ in $U$ which is a contradiction
because we have assumed $\overline{\Sigma}_{\lambda} \cap \overline{O}_{\infty}=\emptyset.$

Now assume $n=2$. Since $\overline{\Sigma}_\lambda \cap \partial \Omega =\emptyset$ and $O_0$ has only one connected component, there exists two distinct point $a,b \in \bar{\Sigma}_\lambda \cap \partial O_0$ such that
\[\partial O_0 \setminus \{a,b\}= V_1 \cup V_2.\]
Now notice that $\Sigma_\lambda \cup V_1$ is a continuous closed curve in $\R^2$. By the Jordan Curve Theorem there exists a bounded open set $U_1$ such that $\partial U_1=\Sigma_\lambda \cup V_1$. Define $U=U_1 \setminus \bar{O}_0\neq \emptyset$.  Then $\partial U \subset \Sigma_\lambda \cup \partial O_0$ which is a contradiction in view of (\ref{cont1}).

In both cases (I) and (II) the contradiction follows from the assumption that
$\overline{\Sigma}_\lambda \cap
\partial \Omega= \emptyset$.   \hfill $\Box$ \\

Since $u \in C^1(\bar{\Omega}\setminus O_0)$, $u$ can be extended to a function in $C^1(\R^n \setminus O_0) \cap BV(\R^n)$. We will denote the restriction of $u$ to $\Omega^c$ by $f$, again. Let $\overline{u}_1$ be the continuous extension of $u_1$ to $\R^n$ with $\bar{u}_1=f$ on $\Omega^c$ and define
\[F_\lambda=\{x\in\R^n\setminus \bar{O}_0: \ \ \bar{u}_1(x)\geq \lambda\}.\]
\begin{remark}\label{openSuperLevelSets}
Let $\Lambda \subset range(u_1)$ be the set defined by Lemma \ref{sameLevelSets} and $\lambda \in \Lambda$.  By Lemma \ref{sameLevelSets} every connected component of $\partial F'_\lambda \cap (\Omega \setminus Z)$ is a $C^1$ hypersurface, where $F'_\lambda$ is defined by $(\ref{Eprime})$. Therefore without loss of generality we may assume that $F_\lambda \cap (\Omega \setminus Z)$ is open, since otherwise $F'_\lambda \cap (\Omega \setminus Z)$ can be replaced by $int (F'_\lambda) \cap (\Omega \setminus Z)$ which differs from $F'_\lambda \cap (\Omega \setminus Z) $, and hence $F_\lambda \cap (\Omega \setminus Z)$ on a set of measure zero. 
\end{remark}

The proof of the following lemma is very similar to that of Theorem 3.7 in \cite{sternberg_ziemer92}. We include the proof for the convenience of the reader. 
\begin{lemma}\label {SZ-lemma}
Let $\Omega$ be a bounded domain with connected Lipschitz boundary. If $x \in \partial ^*F_\lambda \cap \partial \Omega $ and
\[\lim _{r\rightarrow 0}\dashint _{B_r(x)\cap \Omega}|\bar{u}_1(y)-f(x)|dy=0, \]
then $\lambda=f(x)$.
\end{lemma}
{\bf Proof:} Assume $f(x)<\lambda$. Then 
\begin{eqnarray*}
0&=& \lim_{r\rightarrow 0} \frac{1}{|B_r(x)\cap \Omega|} \left( \int_{B_r(x)\cap \Omega\cap \{\bar{u}_1< \lambda\}} |\bar{u}_1(y)-f(x)|dy+ \int_{B_r(x)\cap \Omega\cap \{\bar{u}_1 \geq \lambda\}} |\bar{u}_1(y)-f(x)|dy\right)\\
&\geq& \limsup_{r\rightarrow 0}\frac{1}{|B_r(x)\cap \Omega|}  \int_{B_r(x)\cap \Omega \cap  \{\bar{u}_1\geq \lambda\}} |\bar{u}_1(y)-f(x)|dy\\
&\geq& (\lambda -f(x)) \limsup_{r\rightarrow 0} \frac{|B_r(x)\cap \Omega \cap \{\bar{u}_1 \geq \lambda\}|}{|B_r(x)\cap \Omega|}.
\end{eqnarray*}
Consequently 
\[\limsup_{r\rightarrow 0} \frac{|B_r(x)\cap \Omega \cap \{u_1 \geq \lambda\}|}{|B_r(x)\cap \Omega|}=0.\]
On the other hand since $f$ is the trace of $\bar{u}_1\in BV(\R^n \setminus \Omega)$ on $\partial \Omega$, with a similar argument we conclude that
\[\limsup_{r\rightarrow 0} \frac{|B_r(x)\cap (\R^n \setminus \Omega)\cap {\{\bar{u}_1}\geq \lambda\}|}{|B_r(x)\cap (\R^n\setminus \Omega)|}=0.\]
Therefore 
\[\lim_{r\rightarrow 0} \frac{|B_r(x) \cap \{\bar{u}_1 \geq \lambda\}|}{|B_r|}=0\]
and hence $x \not\in \partial ^* E_\lambda$ which is a contradiction. Similarly $f(x)>\lambda$ leads to a contradiction. Thus $f(x)=\lambda$. \hfill $\Box$ \\

\begin{proposition} \label{takeTheRightValueOnBoundary}
Let $(f,a)$ be an admissible pair and $u_1$ be a minimizer of $(\ref{mainProb})$.  Then for  almost every $\lambda \in \Lambda$
 \[u(\partial F_{\lambda}\cap (\overline{\Omega}\setminus Z))=\{\lambda\},\]
 where $\Lambda$ and $Z$ are defined by Lemma \ref{sameLevelSets} and (\ref{Z}), respectively. 
\end{proposition}
{\bf Proof:} In view of Remark \ref{openSuperLevelSets}, we may assume that $F_\lambda \cap (\Omega\setminus Z)$ is open and every connected component of $\partial F_\lambda \cap (\Omega \setminus Z)$ is a $C^1$ hypersurface intersecting  $\partial \Omega$. Now let $\Sigma$ be a connected component of $\partial F_\lambda \cap (\Omega \setminus Z)$. By Proposition \ref{intersectBoundary}, $\bar{\Sigma}\cap \partial \Omega \neq \emptyset$. Let $x_0\in \bar{\Sigma}\cap \partial \Omega \neq \emptyset$. Since $x_0 \not \in Z$, $|\nabla u(x_0)|> 0$. Hence $x \in \partial^* F_\lambda \cap \partial \Omega$, and by Lemma \ref{SZ-lemma} and Proposition \ref{sameLevelSets} we conclude that $u(\Sigma)=\{\lambda\}$. \hfill $\Box$
\\

It is now straightforward to deduce uniqueness from the results established above. To make the argument rigorous it helps to work with super level sets of the solutions as in \cite{JMN} and \cite{sternberg_ziemer92}. Note however that we do not rely on maximum principles for minimum surfaces that are at the core of the proofs in \cite{JMN} and \cite{sternberg_ziemer92}, but rather on  Lemma \ref{sameLevelSets} and Proposition \ref{intersectBoundary}.

\medskip

{\bf Proof of Theorem \ref{mainResult}:} First we prove that $u_1=u$ a.e. in $\Omega \setminus (Z\cup \bar{O}_0)$. Suppose this is not true, then without loss of generality  we may assume that  there exists $\alpha>0$ such that
\[\mathcal {H}^{n}(N)>0,\]
where
\[N:=\{x\in \Omega\setminus (Z\cup \bar{O}_0): \ \ u_1(x)\geq  u(x)+\alpha\},\]
because otherwise $(f,a)$ can be replaced by the admissible pair $(-f, a)$.
Let
\[
\lambda ^* =\sup \{\lambda: \ \ \mathcal {H}^{n} (\{x\in \Omega\setminus (Z\cup \bar{O}_0): \ \ u(x)\geq \lambda\} \cap N )\geq \frac{\mathcal {H}^{n}(N)}{2}\},
\]
Since $u\in L^1(\Omega \setminus \bar{O}_0)$, $\lambda^* <\infty$.  Define
\[E_1= \{x\in \Omega\setminus  (Z\cup \bar{O}_0) : \ \  u_1(x)\geq \lambda^*+ (1-\beta) \alpha\},\]
then  by Lemma \ref{sameLevelSets} and Proposition \ref{intersectBoundary} there exists $0<\beta<1$ such that $\lambda^*+(1-\beta)\alpha \in \Lambda$.  Also it follows from the definition of $\lambda^*$ that $\mathcal {H}^{n}(K)>0$, where
\[K:=\{x\in \Omega\setminus  (Z\cup \bar{O}_0): \ \ \lambda^*-\beta \alpha < u(x) < \lambda^* \} \cap N.\]

Now let $E_2=\{x\in \Omega\setminus (Z\cup \bar{O}_0): \ \ u(x)\geq \lambda^* \}$.  It is easy to see that $K \subset E_1 \setminus \bar{E}_2 \subset \Omega \setminus (Z\cup \bar{O}_0)$. On the other hand by Remark \ref{openSuperLevelSets} we may assume that $E_1$ is open and hence $E_1 \setminus \bar{E}_2$ is a non-empty open set. Also 
\[\partial (E_1\setminus \bar{E}_2) \subset( \partial E_1 \cap \overline{E_2^c}) \cup (E_1 \cap \partial E_2)\]
and in particular, $\partial (E_1\setminus \bar{E}_2) \subset \partial E_1 \cup \partial E_2$. Notice that $\partial (E_1\setminus \bar{E}_2) \not \subset \partial E_2$, because otherwise $u=\lambda^*$ in $E_1\setminus \bar{E}_2$ which is in contradiction with the assumption $E_1\setminus \bar{E}_2 \subset (\Omega\setminus  Z)$. Let
\[x_0 \in \partial (E_1\setminus \bar{E}_2) \setminus \partial E_2. \]
Then $x_0 \in \partial E_1 \cap \overline{E_2^c}$. By Lemma \ref{sameLevelSets} and Proposition \ref{intersectBoundary} we have
\begin{equation}\label{last}
u(x_0)\in u(\partial E_1)=\{ \lambda^*+(1-\beta)\alpha\}.
\end{equation}
On the other hand
\[u(x_0) \in u(\overline{E_2^c}) \subset (-\infty, \lambda^*]\]
which is in contradiction with (\ref{last}). Hence $u_1=u$ a.e. in $\Omega \setminus (Z\cup \bar{O}_0)$. 

Now let $\Sigma$ be a connected component of $Z$. By the admissibility assumption, $int (\overline{u(O_\infty)})=\emptyset$ and $u$ is continuous.  Therefore $u$ must be constant on $\Sigma$. Since $u=u_1$ in $\Omega \setminus (Z\cup \bar{O}_0)$ and $u_1$ minimizes (\ref{mainProb}), $u=u_1$ a.e. in $\Sigma$. The proof is now complete. \hfill $\Box$

\section{Acknowledgements}

We would like to thank Robert L. Jerrard from whom we have learned a lot through the course of this project.  This work originated and was done in part during the third author's participation in the semester long Thematic Program on Inverse Problems and Imaging in the Fields Institute, January-May, 2012; he was also supported by the sabbatical program at the University of Central Florida, Orlando,  USA,  which made this visit possible. The paper was completed during the second authors participation in the program on Inverse Problems and Applications at the Mittag-Leffler Institute. We would like to thank all three organizations for their support. The first author was partially supported by MITACS and NSERC postdoctoral fellowships. The second author was partially supported by an NSERC Discovery Grant.

 \end{document}